\documentclass[a4paper,10pt]{article}

\usepackage{amssymb}
\newcommand{\inc}{\ensuremath{\subset}}
\newcommand{\demi}{\ensuremath{\frac{1}{2}}}
\newcommand{\gene}{\ensuremath{\mu : M\to \mathcal{C}^*\times (S^1)^r}}
\newcommand{\z}{\ensuremath{\mathbb{Z}}}
\newcommand{\T}{\ensuremath{\mathcal{T}}}
\newcommand{\orb}{\ensuremath{\mathcal{O}}}
\newcommand{\rev}{\ensuremath{\mathcal{V}}}
\newcommand{\esp}{\ensuremath{\mbox{\hspace{2cm}}}}
\newcommand{\G}{\ensuremath{\mathcal{G}}}
\newcommand{\C}{\ensuremath{\mathcal{C}}}
\newcommand{\Q}{\ensuremath{\mathbb{Q}}}
\newcommand{\R}{\ensuremath{\mathbb{R}}}
\newcommand{\cqfd}{\ensuremath{\bigcirc}}
\newcommand{\fin}{\begin{flushright}
\cqfd
\end{flushright}}
\newcommand{\matr}[2]{\ensuremath{\left( \begin{array}{#1}#2 \end{array}\right)}}
\newcommand{\sm}{\ensuremath{(M,\omega)}}
\newcommand{\om}{\ensuremath{\omega}}
\newcommand{\smp}{\ensuremath{(M,\omega')}}
\newcommand{\sympo}{\ensuremath{\mbox{Symp}_o\sm}}
\newcommand{\ham}{\ensuremath{\mbox{Ham}\sm}}
\newtheorem{theo}{Theorem}[section]
\newtheorem{exa}{Example}[section]
\newtheorem{cor}{Corollary}[section]
\newtheorem{prop}{Proposition}[section]
\newtheorem{lem}{Lemma}[section]

\newtheorem{rem}{Remark}[section]

\newtheorem{defi}{Definition}[section]
\newcommand{\ddt}{\ensuremath{\frac{\mathrm{d}}{\mathrm{dt}} \left|_0 \right. }}
\begin{document}
\title{On Generalized Moment Maps for Symplectic Compact Group Actions }
\author{Pierre Sleewaegen\thanks{in collaboration with the Math. Department at U.L.B. Brussels}\\
 165, Av. Winston Churchill\\
B-1180 Brussels - BELGIUM\\
p.sleewaegen@tiscali.be      }
\date{April 2003}
\maketitle

\begin{abstract}
In this paper a moment map is defined for an arbitrary symplectic action of a compact connected Lie group
on a closed symplectic manifold, in the spirit of the circle-valued map introduced by
McDuff in the case of non-Hamiltonian circle actions. We show that, for torus actions, the 
Atiyah-Guillemin-Sternberg convexity theorem remains valid in our context. Also, 
we study the equivariance properties of generalised moment map and some facts about
Marsden-Weinstein reduction procedure, allowing e.g. to reformulate a proof
of Kim's result that "complexity one" symplectic torus actions are Hamiltonian.
As illustration of the use McDuff moment maps, we give a symplectic proof of the finiteness
of certain symmetry groups of genus $\geq 2$\ compact oriented surfaces.
\end{abstract}

Generalized moment maps for non-Hamiltonian symplectic $S^1$-actions on closed symplectic manifolds were introduced by McDuff in
\cite{DuMMCA1988}, allowing her to prove that, for closed symplectic 4-manifolds,
any $S^1$-action with a fixed point is Hamiltonian. A McDuff moment map is defined as follows :
\begin{defi}
Let \sm be a closed symplectic manifold endowed with a symplectic non-Hamiltonian $S^1-$action.
Let $\eta\in s_1\cong\mbox{\it Lie}(S^1)$\ generate $\ker\exp\cong\mathbb{Z}$,
$X_\eta$\ be the corresponding fundamental field on \sm , 
$\alpha$\ be the basis of $s_1^{*}$\ dual to $\{\eta\}$\ generating the invariant 1-form
$\tilde\alpha$\ on $S^1$\ . Then a McDuff moment map is a map
$$\mu : M\to S^1\mbox{\hspace{2cm} with\hspace{2cm}}
\mu^*\tilde{\alpha}= i_{X_{\eta}}\om'$$
\end{defi}
where $\om'$\ is an integral symplectic form constructed from \om .\\ (more details will be given in
section 1 below)

Important properties of McDuff moment maps are (see \cite{AuToTo} p.77) :
\begin{itemize}
\item $\mu$\ is $S^1$-invariant
\item $\mu(M)=S^1$\ (surjectivity)
\item each fiber of $\mu$\ has the same number of connected components
\item $\mu$\ admits no local extremum 
\end{itemize}

The motivation behind the present paper is to gain a basic understanding of the possible " moment map "
obtained when McDuff's construction is extended to compact group actions
(extensions to non-compact groups belong to an other story).
Here is an outline of the work :\\

We start with an application of the $S^1$\ case : using a result of Ono \cite{OCGA} based on
McDuff moment map, we will obtain a symplectic proof of the following well-known result :
\begin{theo}
Let S be a smooth compact oriented genus $\geq 2$\ surface. Then
\begin{enumerate}
\item for any Riemann metric on S, the isometry group is finite.
\item for any symplectic structure \om\ on S, the symplectomorphism group of $(S,\om)$ has no connected compact Lie subgroup.
\item for any complex structure on S, the group of holomorphic transformations is finite.
\end{enumerate}
\end{theo} 

As generalization of McDuff moments, we introduce for any symplectic action $\sigma$\ of
a connected compact Lie group $G$\ on a compact symplectic manifold \sm\ a
 "generalized moment"\footnote{It is worth to draw a parallel between the present work and the two papers \cite{GiCMM} and \cite{GiSACG}.}
 .
Its existence relies on the facts that
\begin{itemize}
\item  (see section 2) up to a covering $G$\ decomposes as a direct product
 $G=C\times T^r$\ where $C$\ is the maximal connected compact subgroup of $G$\
on which $\sigma$\ is Hamiltonian, and $T^r$\ is a torus in $G$\ acting in a 
totally non-Hamiltonian way
\item (section 3) an invariant integral symplectic form $\om'$\ exists on $M$\ , relative to which
$\sigma$\ is exactly as Hamiltonian as it is relative to \om\ .
\end{itemize}
legitimating the
\begin{defi}
A generalized moment map for $\sigma$\ is a map
$$ \mu : M\to \mathcal{C}^*\times (S^1)^r $$
where the first factor $\mu_1 : M\to \mathcal{C}^*$\ is a genuine Hamiltonian
moment map and the r other factors are McDuff moments associated to each circle
factor of a decomposition $T^r\cong (S^1)^r$\ , all this evaluated relative to
a unique integral invariant symplectic form $\om'$\ 
\end{defi}

There natural notion of equivariance for generalized moment maps : $\mu$\ is said to be equivariant
if $\mu_1$\ is $Ad^*(C)$-equivariant and $T^r$-invariant, while $\mu_2:M\to (S^1)^r$\ is $G$-invariant.
For an equivariant generalized moment, the Marsden-Weinstein reduction procedure
can be xeroxed from the Hamiltonian case.  But the basic example of the 2-torus acting on itself
with moment $\mu:(\mathbb{R}/\mathbb{Z})^2\to(\mathbb{R}/\mathbb{Z})^2:[p,q]\to[-q,p]$\
shows that,  while equivariance can always be assumed in the Hamiltonian case (\cite{GoTuSMP1991} ,prop.2),
it fails here.  However, we will show that things are nicer in the presence of fixed points:
\begin{prop}
If any fixed point exists for the $T^r$-action on $M$\ , $\mu$\ can be chosen equivariant .
\end{prop}
This is not so surprising : in many (but not every !) situations that we recall in section 3, the presence of a fixed point
even implies that the action is Hamiltonian.\\
When the moment map is not equivariant in the above "coadjoint" sense, we will
show using arguments as in \cite{LiMaSM} and \cite{CoDaMo1988} that some "affine" action $Aff$\
exists on $(S^1)^r$\ regarding to which $\mu_2$\ is equivariant. This also authorizes reduction procedure, 
in the less usual way described in \cite{LiMaSM} .

Section 4 is devoted to a non-Hamiltonian version of the Atiyah-Guillemin-Sternberg convexity theorem (\cite{AtCCH1982,GuStCPI1982}):
\begin{theo}
When $G$\ is a torus, $M$\ compact connected, let $\Delta=\mu_1(M)$\ be the 
convex polytope image of $\mu_1$\ . Then
$$\mu(M)=\Delta \times (S^1)^r$$
Moreover $\mu :M\to \mu(M)$\ is an open map.
\end{theo}
We have the 
\begin{cor}
Let $b_1(M)$\ denote the first Betti number of $M$. Then
$$ r\leq b_1(M) $$ 
\end{cor}
Also, we deduce from the convexity theorem some "heredity" property : when $T_1\times T_2$\ acts on \sm, $T_1,T_2$\ being torus groups, with a naturally equivariant moment,
and  the $T_2$-action is non-Hamiltonian, the induced $T_2$-action on a $T_1$-reduced space
$M_{red}^{T_1}$\ remains non-Hamiltonian. This allows us to reformulate in the Appendix a proof
of a result of M.K. Kim generalizing the introducing McDuff result to "complexity one" symplectic action, i.e.
the situation where a n-torus acts symplectically on a (2n+2)-dimensional \sm.\\

{\bfseries{ Notation conventions}}: throughout this paper, \sm\ will denote a 
closed connected symplectic manifold, $G,C$\ compact connected Lie groups (as far
as generalized moments are around), $T$\ a torus group,
$\mathcal{G},\mathcal{C},\mathcal{T}$\ their respective Lie algebras, $\mathcal{G}^*,\mathcal{C}^*,\mathcal{T}^*$\
their duals.  To each element $\xi$\ in the Lie algebra of a group acting on $M$\
$X_{\xi}$\ will be the associated fundamental field on $M$\ : $X_{\xi}(x)=\ddt (\exp -t\xi .x)$\ .
Finally, by $>S<\subset V$\ we denote the vector space generated by a subset $S$\ in a vecor space $V$.

\section{McDuff moments and genus $\geq 2$\ surfaces}
Let $\sigma:S^1\times M\to M$\ be a symplectic non-Hamiltonian $S^1$-action on \sm\ .
The construction of a McDuff moment map goes as follows : an invariant
symplectic form $\om''$\ with a rational cohomology class is chosen near \om\ ,
preserving the non-Hamiltonian property.  Then $\om'=k.\om''$\  is in an integral class
for some $k\in\mathbb{Z}$\ . Now, trying to define the moment map in the usual manner as
$$\mu^{usu}(x)=\int_{\gamma_x} i_{X_{\eta}}\om' \mbox{\hspace{1cm}}x\in M, \gamma_x:[0,1]\to M,
\gamma_x(0)=\tilde{x},\gamma_x(1)=x$$
for some base point $\tilde{x}$\ , we see that the obstructions lie in the non-vanishing
 of the periods $\int_C i_{X_{\eta}}\om'$\ for some 1-cycles in M.  But those periods are
 in \z\ , because $C$\ can be pushed to a 2-cycle $\sigma(C):[0,1]^2\to M:(s,t)\to
 (\exp s\eta).C(t)$\ and
 $$ \z\ni \int_{\sigma(C)}\om'=\int_{[0,1]^2}\om'_{(\exp s\eta).C(t)}(X_{\eta}({(\exp s\eta).C(t)}),(\exp s\eta)_*(\dot{C}(t))=
 \int_C i_{X_{\eta}}\om'$$
where the invariance and the integrality of $\om'$\ was used.  The McDuff moment $\mu$\ is then simply
$\mu^{usu}$\ seen as a map to $\R/\z$\ .
The properties of $\mu$\ recalled in the introduction are either immediate or consequence 
of the fact that $\mu$\ is a Morse-Bott function with critical points of even indices (see \cite{AuToTo}).
\begin{rem}[Aside]
It is possible to render $\mu$\ fiber-connected.
\end{rem}
Indeed, the non-existence of a local extremum for $\mu$\ implies that $\mu$\ is an open map
onto $S^1$\ , and defining as in \cite{CoDaMo1988,HiNePl1994} the equivalence 
relation $\sim$\ on $M$\ by saying $x\sim y$\ if $\mu(x)=\mu(y)$\ and $x$\ and $y$\ belong to 
the same connected component of $\mu^{-1}(\mu(x))$\ , we see that $\tilde{\mu}:(M/\sim)=\tilde{M}\to S^1$\ is a finite covering.
Such a covering is always of the form $\pi:S^1 \to S^1 :t\to t^k$\ for some $k\in\z$, whence $\mu$\ factorizes as
$$M\stackrel{\mu^{fc}}{\longrightarrow} \tilde{M}\cong S^1 \stackrel{\pi}{\longrightarrow} S^1$$
where $\mu^{fc}$\ is fiber-connected.  We will discuss implications of this observation below in section 4 and in the
Appendix.

As example of the use of McDuff moment maps, we will now provide a proof of Theorem 0.1.
Compact connected group actions on surfaces are naturally not very mysterious (see e.g. the full
classification of $S^1$-actions on compact surfaces in \cite{AuToTo}), the point here is to show 
that a symplectic argument can be relevant.  The starting point is the following
theorem of Ono, obtained from the properties of McDuff moments :
\begin{theo}[\cite{OCGA}]
Let \sm\ be a closed symplectic manifold.
\begin{enumerate}
\item if the second homotopy group $\pi_2(M)$\ vanishes, then a circle group action
on $M$\ preserving \om\ has no fixed point.
\item Moreover if every abelian subgroup of $\pi_1(M)$\ is cyclic, there is no circle group
action on $M$\ preserving \om\ .

\end{enumerate}
\end{theo}
It is also useful to collect here some standard facts about (symplectic) $G$-actions on \sm\ :
\begin{itemize}
\item isotropy groups $G_x=\{g\in G|g.x=x\}$\ are divided into conjugation classes characterizing the
"orbit type" of point in $M$\ . 
\item There is only a finite number of orbit types (\cite{GuStSy}, prop.27.4 or \cite{AuToTo} prop. 2.2.3)
\item There exists a "principal" orbit type s.t. the corresponding points (the "principal stratum" $M_{princ}$\ ) in $M$\ form a dense open subset.
\item For $H$\ a closed subgroup of $G$, define
$$ M_H=\{x\in M|G_x=H\}.$$
Then each connected component of $M_H$\ is a symplectic submanifold of $M$\ (\cite{GuStSy} p.203)
\end{itemize}

{\bfseries{Proof of theorem 0.1 :}} We remark that contradicting the conclusion of
the theorem amounts to the existence of some non-trivial $S^1$-action $\sigma$\ preserving a symplectic form
\om\ on $S$\ : for 1. it follows from the fact that the isometry group of a Riemann structure
on $S$\ is a compact Lie transformation group (\cite{KoTrGr}, p.39) preserving the induced (symplectic) volume form ;
for 2. it is obvious ; for 3. we get an holomorphic transformation invariant metric $g$ by
use of the uniformization theorem (\cite{JoRiSu}, theorem 4.4.1) and the fact that holomorphic transformations
of the hyperbolic half-plane are isometries (\cite{JoRiSu},pp.28-29). \om\ is again the Riemannian volume $v_g$\ .

Since $\pi_2(S)=0$\ , part 1. of the Ono theorem applies and $\sigma$\ has no fixed point.
From the fact that the $M_H$\ are symplectic and in particular even dimensional, we deduce that
$S=S_{princ}$\ , i.e. all points have the same isotropy group $\cong \z/k\z$\ for some $k$ . $\sigma$\ can then be replaced by a free action,
so that $S$\ becomes a principal $S^1$-bundle over a 1-dimensional compact base, forcing
$S\cong S^1\times S^1$\ contradicting the genus hypothesis.
\begin{flushright}
\cqfd
\end{flushright}
\begin{rem}
Hurwitz's theorem (\cite{KoTrGr} p.88) states that the order of the group of holomorphic transformations of $S$\ is at most
$84(g-1)$\ , and it is even known that it reduces to identity when $g\geq 3$\ (\cite{GrHaPr} p.276).  Naturally, our
technique doesn't give access to those finer results.
\end{rem}
\section{The basic lemmas}
The two lemmas below will legitimate the definition of generalized moments
in the next section.
\begin{lem}
Let $H$\ be a dense connected Lie subgroup of G, and suppose
given an Hamiltonian action of $H$\ on \sm\ . If this action extends to $G$, 
the extension is Hamiltonian too.
\end{lem}
\begin{rem}
No question of equivariance is considered until the next section.
\end{rem}
The tools we need are :
\begin{enumerate}
\item Symplectic actions are always Hamiltonian when the group is semisimple (by
the Whitehead lemmas and Theo. 26.1 in \cite{GuStSy})
\item Up to a covering that is not relevant here, $G=K\times T$\ , where $K$\ is 
compact semisimple and $T$\ is a torus (\cite{ZeCoLi} p.297)
\item If a torus $T$\ acts symplectically on \sm\ and $S^1 \subset T$\ is a circle
subgroup acting in a non-Hamiltonian way, then the modified symplectic form $\om'$\
used to define the associated McDuff moment can be chosen $T$-invariant. (this follows
from the starting argument of the proof of the next lemma).
\item Near each fixed point $x$\ of a symplectic $T$-action, there exists a local moment
map which reads in an appropriate Darboux chart $(x_1,\ldots,x_n,y_1,\ldots,y_n)$\ as
$$\mu(x_i,y_j)=\mu(0) + \demi\sum_{i=1}^n \alpha_i (x_i^2+y_i^2)$$
where the $\alpha_i\in\mathcal{T}^*$\ are the weights of the linearized $T$-action on
the tangent space $T_xM$\ (\cite{GuStSy} p.250). Those weights are independant of the choice of \om\ or $\om'$\ above
and belong to the weight lattice in $\mathcal{T}^*$, meaning that they take integral values on the
$\z^n$\ lattice $\ker\exp$\ in $\mathcal{T}$\ .
\end{enumerate}
{\bfseries{Proof of lemma 2.1}}:
By continuity considerations, the $G$-action is symplectic. Also, by tools 1) and 2)
above, it is readily seen that the question reduces to the case of a dense 1-parameter
subgroup $H$\ in a torus $T$, the $H$-action being Hamiltonian. The associated Hamiltonian
$f_H:M\to\R$\ takes a minimum value ($M$\ is compact) at some $H$-fixed point 
$p\in M$\ . Then $p$\ is $T$-fixed, and in the notations of tool 4. above $f_H$\ can be written near $p$\ as
$$f_H=f_H(p)+ \demi\sum_{i=1}^n \alpha_i(h) (x_i^2+y_i^2)$$
with $\alpha_i(h)\geq 0\forall i$\ , $h$\ being some generator of the Lie algebra of $H$. Since
the coordinates of $h$\ are $\Q$-linearly independent in a natural basis $\{e_1,\ldots,e_r\}$\ of 
\T\ on which the $\alpha_i$\ take integral values, we see that $\alpha_i(h)=0$\ is only possible
if $\alpha_i=0$\ in $\T^*$\ , and that $\alpha_i(h)>0$\ otherwise. Now we can choose
linearly independent vectors $r_1,\ldots,r_r$\ in \T\ with rational coordinates in the
basis $\{e_k\}$\ and close enough to $h$\ to have $\alpha_i(r_j)>0\;\;\forall j,\alpha_i\neq 0$.
Each $r_j$\ generates a circle subgroup $C_j$\ of $T$.
The lemma will be proved if we show that every $C_j$-action on $M$ is Hamiltonian.
If not then the action of at least one circle, say $C_1$, admits a McDuff moment
$\mu_1:M\to S^1$. But, using the independence of the weights to a small change
of symplectic structure in tool 4. above, we see that near $p$\ and in well-chosen
coordinates on $M$\ and $S^1$, $\mu_1$\ looks like
$$\mu_1(p)+\demi\sum_{i=1}^n \alpha_i(r_1) (x_i^2+y_i^2)$$
furnishing a local extremum of $\mu_1$. This is impossible.
\begin{flushright}
\cqfd
\end{flushright}
\begin{rem}
This lemma is perhaps true even without compactness assumption on $G$.
\end{rem}
In fact it is linked with the important Flux Conjecture (see \cite{LMPFC}) that
we now restate :
\begin{defi}
an Hamiltonian isotopy on \sm\ is a family $\Phi_t$\ of symplectomorphisms
obtained from a smooth family of Hamiltonian functions $H_t:M\to\R\;\;t\in [0,1]$\
as the flow of the time-dependent vector field $X_t$\ s.t. $i_{X_t}\om=dH_t$.
\end{defi}
The set $\ham=\{\Phi_1|\Phi_t \mbox{ is a Hamiltonian isotopy }\}$\ is a normal path-connected
subgroup of the identity component \sympo\ of the symplectomorphism group (see \cite{DuSaIn} p.311).
The Flux Conjecture asserts :\\

" \ham\ is $C^1$-(or/and)$C^0$-closed in \sympo\ "\\

Our lemma for more general $G$\ would be a corollary of the conjecture as follows :
since each element in $H$\ belongs to \ham\, the Flux Conjecture gives that the entire
$G$\ belongs to \ham. Now every smooth path $\Psi_t\in\ham$\ is generated by Hamiltonian vector fields
(Prop. 10.17 in \cite{DuSaIn}). This implies that $\forall X\in\mathcal{G}$, the 1-parameter
subgroup $\{\exp tX\}$\ in $G$\ corresponds to a Hamiltonian isotopy $\Phi_t^X$\ s.t.
$\Phi^X_{t+t'}=\Phi^X_t\circ\Phi^X_{t'}$\ . Then it is easy to show that the corresponding Hamiltonian $H_t^X$\
is time-independent, providing the desired Hamiltonian function $H^X$ .\\

Back to the compact $G$\ case and a symplectic action $\sigma$\ on \sm\ remember that every subtorus $T^i$\ in a torus $T^n$\ has a complementary subtorus
$T^{n-i}$\ with $T^n=T^i\times T^{n-i}$. Then we deduce from lemma 2.1  that (up to an irrelevant covering) $G$\ decomposes as a direct product $C\times T^r$\ 
where $C$\ is compact connected, $T^r$\ a r-torus, $\sigma$\ being Hamiltonian on
$C$\ and totally non-Hamiltonian on $T^r$\ (i.e. $i_{X_\eta}\om$\ is non-exact $\forall\eta\in\T$ ).\\
In order to define the generalized moment map, we first want to replace \om\ with
an integral symplectic 2-form s.t. the just mentioned decomposition remains valid.

\begin{lem}
The situation being as just described, there exists an integral $G$-invariant
symplectic form $\om'$\ on $M$\ with $i_{X_\xi}\om$\ exact $\Leftrightarrow$\ $i_{X_\xi}\om'$\ exact,
$\forall \xi\in\mathcal{G}$.

\end{lem}
{\bfseries{proof :}} borrowing an argument from \cite{GoTuSMP1991}, we deduce from the compacity
of $M$\ that it is of "finite integral rank" i.e. 
$$H^i(M,\R)\cong H^i(M,\z)\otimes\R$$
as a consequence of the universal coefficients formula in cohomology when the
cochain complex is finitely generated (Ex.2 p.172 in \cite{McLaHo}).\\
So we can write \om\ as
$$\om=a_1\om_1+\ldots+a_l\om_l\mbox{\hspace{2cm}}a_i\in\R,\om_i\in H^2(M,\z)$$
and since averaging over $G$\ has no consequence on the cohomology class, we 
can furthermore assume that each $\om_i$\ is $G$-invariant. Now let $\{\xi_1,\ldots,\xi_c\}$\
be a basis of the Lie algebra \C\ of $C$\ , consisting of circle generators satisfying
$\exp_G(\xi_i)=e,\exp_G(t\xi_i)\neq e \forall t\in (0,1)$\ . Let $\{\eta_1,\ldots,\eta_r\}$\
be such a basis for \T, and denote $\{\mu_k\}$\ the union of these two basis.
Let $\gamma_1,\ldots,\gamma_s$\ be smooth 1-cycles generating the 1-homology of $M$.
Then the Hamiltonian character of the action is reflected in the properties
$$(A)\mbox{\hspace{2cm}}\int_{\gamma_i}i_{X_{\xi_j}}\om=0=\sum_{k=1}^la_k\int_{\gamma_i}i_{X_{\xi_j}}\om_k$$
$\forall i=1\ldots s,\forall j=1\ldots c$\ and also $\forall j=1\ldots r,\exists k_j\in\{1\ldots s\}$\ s.t.
$$(B)\mbox{\hspace{2cm}}\int_{\gamma_{k_j}}i_{X_{\eta_j}}\om\neq 0$$
Now we can as in section 1 push the 1-cycles with the various circle actions to obtain
various 2-cycles $\sigma_{\mu_j}\gamma_i$\ and prove that 
$$\forall i,j,k,\mbox{\hspace{2cm}}\int_{\gamma_i}i_{X_{\mu_j}}\om_k \in\z$$
So $(a_1,\ldots,a_l)^t\in\R^l$\ is a non-zero solution of a system of linear equations
of the form $A.a=0\;\;\;,\;\;\;A\in\z^{s\times l}$. The solution space of this system
can thus be generated by vectors in $\Q^l$.  We may then replace $(a_1,\ldots,a_l)$\
by an arbitrary close rational vector $(q_1,\ldots,q_l)$\ s.t.
$$\om''=q_1\om_1+\ldots+q_l\om_l$$
remains non-degenerate and still satisfies properties (A) and (B).
Choosing $\om'=k.\om''$\ for a large enough integer $k$\ concludes this proof.
\begin{flushright}
\cqfd
\end{flushright}

\section{Definition and equivariance of the generalized moment}
Everything to guarantee the existence of generalized moment maps is now established.
Keeping the notations of the preceding section, here is again the definition :
\begin{defi}
A generalized moment map for $\sigma$\ is a map
$$ \mu : M\to \mathcal{C}^*\times (S^1)^r $$
where the first factor $\mu_1 : M\to \mathcal{C}^*$\ is a genuine Hamiltonian
moment map and the r other factors are McDuff moments associated to each circle
factor of a decomposition $T^r\cong (S^1)^r$\ , all this evaluated relative to
a unique integral invariant symplectic form $\om'$\ 
\end{defi}
This definition is the one that appears to us the closest in spirit to the original McDuff's construction.
\begin{rem}
Many choices are implicit in this definition, but we will study properties unaffected by those choices.
\end{rem}
\begin{rem}
The possibility of defining a moment map taking values in some cylinder obtained
by quotienting $\G^*$\ adequately is already discussed in \cite{CoDaMo1988}, in
a completely general situation. The problem is to prevent this cylinder from
being reduced to a point ; this is the reason for our integrality obsession here, and
for our seemingly limited framework.
\end{rem}
One of the important constructions allowed in the Hamiltonian case is Marsden-Weinstein
reduction, working as follows (our reference here is \cite{LiMaSM}) :\\
let $\mu : M\to\G^*$\ be the moment map associated with an Hamiltonian action of
a Lie group $G$\ on \sm. Then there exists a (unique) affine action $a_\theta$\ of
$G$\ on $\G^*$\ defined by a 1-cocycle $\theta:G\to\G^*$\ for the coadjoint
action, giving the $\mu$-equivariance ($\mu(g.x)=a_\theta(g)\mu(x)$).\\
Given a weakly regular value $v\in\mu(M)\subset\G^*$\ , $G_v$\ its $a_\theta$\ isotropy group,
$G_v^o$\ its neutral component, the reduced space associated to $v$\ is
$$M_{red}^{v}=\mu^{-1}(v)/G_v^o\mbox{\hspace{1cm}or\hspace{1cm}}
M_{red}^{v}=\mu^{-1}(v)/G_v$$
In many circumstances, $M_{red}$\ is still a manifold (or an orbifold) and has
a naturally induced symplectic structure. This reduction scheme appears as a 
very useful tool in symplectic geometry.\newpage
Now a fundamental observation is that the only global property of the moment used
to construct the reduced space and its symplectic structure is its equivariance.
This leads us to consider equivariance properties of our generalized moments. The starting point
is that, when $G$\ is compact and the action is Hamiltonian, $\mu$\ can always be chosen
$Ad^*$-equivariant, i.e. $\theta=0$\ (\cite{GoTuSMP1991}). This brings our first
equivariance property :
\begin{lem}
Let \gene\ be a generalized moment map for a symplectic $G$-action.\\ Then the first
factor $\mu_1:M\to\C^*$\ can be chosen $Ad^*(C)$-equivariant and $T^r$-invariant, the 
second factor $\mu_2:M\to(S^1)^r$\ being $C$-invariant.
\end{lem}
{\bfseries{proof :}} Starting with an $Ad^*(C)$-equivariant $\mu_1'$, we observe
as in \cite{GoTuSMP1991} that averaging $\mu_1'$\ over the $T^r$-action still
delivers an $Ad^*(C)$-equivariant moment (now $T^r$-invariant!)  $\mu_1$\ for the $C$-action.
The $G$-invariance of $\mu_2$\ is then natural since by the defining property of a
moment map, infinitesimal $T$-invariance of $\mu_1$\ means
$$\om'(X_{\xi},X_\eta)=0\mbox{\hspace{2cm}}\forall\xi\in\C,\forall\eta\in\T$$
But this is equivalent to the infinitesimal $C$-invariance of $\mu_2$.  $C$\ being connected,
this implies the full $C$-invariance of $\mu_2$.
\fin
The equivariance property expected from the Hamiltonian situation is contained
in the following definition
\begin{defi}
$\mu$\ is naturally equivariant if $\mu_1$\ is $Ad^*(C)$-equivariant, $T^r$-invariant,
and $\mu_2$\ is $G$-invariant.
\end{defi}
Natural equivariance is not always obtainable :
\begin{exa}[The 2-torus]
\end{exa} 
Consider the standard action of $T^2\cong\R^2/\z^2$\ on itself :\\
$[p,q].[r,s]=[p+r,q+s]$. The invariant and integral symplectic form\\
$\om="dx\wedge dy/\z^2"$\ gives a generalized moment map
$$\mu:T^2\to T^2:[p,q]\to[q,-p]$$
which can of course not be chosen as $T^2$-invariant.\\

This simple example shows that natural moment equivariance is not the rule for
non-Hamiltonian actions, but there is a situation where things remain nice :
\begin{prop}
If any fixed point exists for the $T^r$-action on $M$\ , $\mu$\ can be chosen naturally equivariant .
\end{prop}
\begin{rem}
The existence of a fixed point for a symplectic torus action simplifies in many cases
the analysis radically by forcing the action to be Hamiltonian.
\end{rem}
This is so e.g. when
\begin{itemize}
\item $M$\ is K\"{a}hler or of Lefschetz type (see \cite{DuSaIn} p.150)
\item $M$\ is 4-dimensional (see \cite{DuMMCA1988})
\item $M$\ is monotone (see \cite{DuSaIn} p.152)
\item The action is semifree with isolated fixed points (see \cite{ToWeSS})
\item The action is effective and $\dim T=\dim M/2$\ (see \cite{GiCMM})
\end{itemize}
Nevertheless the subject of non-Hamiltonian symplectic actions with fixed points
on a compact \sm\ is not empty, as shown by a 6-dimensional example in \cite{DuMMCA1988}.\\

The proof of prop. 3.1 will be based on the following two facts about symplectic actions :
\begin{enumerate}
\item For a torus $T$, any generalized moment $\mu$\ for the $T$-action is $T$-invariant if and only if
the $T$-orbits are isotropic in $(M,\om')$\. (see \cite{AuToTo} prop. 3.5.6)
\item Near any isotropic $G$-orbit \orb\ in \smp\ the equivariant isotropic embedding 
theorem provides a model for a G-invariant neighbourhood of \orb\ in \smp, as described
e.g. in \cite{SjLeSSSR1991}, prop. 2.5 :
\begin{prop}[local normal form for the moment map] 
Let $H$\ be the stabilizer of $p\in\orb$\ and $V$\ be the symplectic slice to the 
orbit \orb. Then a neighbourhood of \orb\ is equivariantly symplectomorphic to
a neighbourhood of the zero section of $Y=G\times_H(m^*\times V)$\ with the equivariant
$G$-moment map $J$\ given by the formula
$$J([g,\mu,v])=Ad^*(g)(\mu+\Phi_V(v))$$
\end{prop}
We will not explain the notations here, the important fact is that any isotropic
orbits possesses an invariant neighbourhood on which the action is Hamiltonian with an $Ad^*$-equivariant
moment.
\end{enumerate}
{\bfseries{Proof of prop. 3.1 :}} By Lemma 3.1, we see that it is enough to prove
$T$-invariance for generalized moment maps of symplectic $T$-actions with fixed points, $T$\ being a torus.
By fact 1) above it is equivalent to show that every $T$-orbit is isotropic.
Clearly, since the isotropy of the orbit through $x\in M$\ is equivalent to
$$\om'(X_\xi,X_\eta)=0\mbox{\hspace{2cm}}\forall \xi,\eta\in\T,$$
the set of points with isotropic orbits in $M$\ is closed. But the Sjamaar-Lerman proposition in fact 2),
applied to the torus case, imply the existence of a (genuine) equivariant moment map
on an invariant neighbourhood of any isotropic orbit, hence again by fact 1) the isotropy of
every orbit in this neighbourhood.  So the set of points with isotropic orbit is
open and closed in \smp, and nonempty by the fixed point hypothesis. $M$\ being connected, the
proposition follows.
\fin

Our next goal is to see that some equivariance in the Libermann-Marle sense above
remains whithout any fixed point assumption. Namely we will prove :
\begin{prop}
There exists a matrix $Z\in \z^{r\times r}$\ with zero diagonal, defining an affine action
$Aff^Z$\ of $T^r$\ on itself, s.t.
$$\mu_2(t.x)=Aff^Z(t)\mu_2(x)\esp \forall t\in T^r, x\in M$$
\end{prop}
\begin{rem}[Definition of $Aff^Z$]

\end{rem}
The affine action $Aff^Z$\ is defined by :
$$ Aff^Z(s_1,\ldots,s_r)(t_1,\ldots,t_r)=(\prod_j(s_j)^{Z_{1j}}.t_1,\ldots,\prod_j(s_j)^{Z_{rj}}.t_r)$$
In the two-torus example, $Z=\matr{cc}{0&1 \\ -1 & 0}$.\\

{\bfseries{proof of prop. 3.3 :}} we see $\mu_2$\ as taking values in $\R^r/\z^r$, and
then providing an ordinary moment $\mu_O:M\to \T^*\cong \R^r$\ on each open set $O\subset M$\ around a point $o$,
s.t. $\mu_2(O)\inc \mu_2(o)+[-\demi,+\demi[^r$.\\
Let us choose the following data : two nested finite covers $\orb=\{\orb_i\}\inc\rev=\{\rev_i\}$\ of $M$\ (compact!),
with $\mu_2|_{\rev_i}$\ an ordinary moment and an inversion-invariant identity neighbourhood
$U\inc T^r$\ small enough to have $t_1.t_2.\orb_i\inc\rev_i\forall t_1,t_2\in U$.
Then a localization of the argument of Th\'eo. 3.2 in \cite{LiMaSM} gives :\\
$\forall i, \forall x \in \orb_i,\forall t\in U, \mu_2|_{\rev_i}(t.x)=
a_i(t,\mu_2|_{\rev_i}(x))$, where 
$a_i(t,v)=Ad^*(t)v+\theta_i(t)$, $\theta_i:U\inc T^r\to \T^*$\ being a 1-cocycle
for the coadjoint action :\\
$\theta_i(t_1.t_2)=\theta_i(t_1)+Ad^*_{t_1}\theta_i(t_2)$.\\
The coadjoint action being trivial in the torus case, we have 
$$a_i(t,v)=v+\theta_i(t)\esp \theta_i(t_1.t_2)=\theta_i(t_1)+\theta_i(t_2)$$
Now, using the connectedness of $M$, we deduce that the $\theta_i$\ in fact are not
dependent on $i$, so there is a unique (local) cocycle $\theta:U\inc T^r\to\T^*$\
associated to $\mu_2$. We have to show that $\theta$\ can be extended from
$U$\ to $T$, keeping its relation to $\mu_2$. Remark that we are interested in the
image of $\theta$\ up to the $\z^r$-lattice. The integrality assumption on $\om'$\ is
crutial at this point : let $\mu_2^1,\ldots,\mu_2^r$\ be the $r$\ McDuff components
of $\mu_2$, $S^1_1,\ldots,S^1_r$\ the corresponding circles in $T^r$,
 $\eta_1,\ldots,\eta_r\in\T$\ the corresponding 'unit' vectors, 
 $\tilde{x}\in M$\ a base point, 
 $C_1=S^1_1.\tilde{x},\ldots,C_r=S^1_r.\tilde{x}$\ r 1-cycles generated by the action.
 Define 
 $$Z_{ij}=\int_{C_j}i_{X_{\eta_i}}\om'\esp\in\z^{r\times r}.$$
 Then it is not hard to see that $\theta$\ induces the homomorphism
 $$T^r\to T^r:s\to Aff^Z(s)(e)$$
 and that $\mu_2$\ is equivariant with respect to $Aff^Z$.
 \fin
\begin{rem}
In the Hamiltonian case, $\theta$\ is global and induces an homomorphism $T^r\to\R^r$, forcing
$\theta$\ to be zero.
\end{rem} 
\begin{rem}
This result was deduced from a corresponding discussion about affine Poisson structures
in \cite{CoDaMo1988}.
\end{rem} 
\begin{cor}
If $Aff^Z$\ is free on some subtorus $T^s\inc T^r$, then $M$\ is a (trivial) $T^s$-principal bundle.
\end{cor}
\begin{cor}
If $Z$\ is of rank r, then the $T^r$-action is locally free, i.e. every point in $M$\ has a finite stabilizer.
\end{cor}

Those equivariance properties in principle give access to the machinery of symplectic
reduction along the Hamiltonian lines. It may seem strange physically to look at 
objects constructed from a different symplectic structure (unless \om\ is integral); 
for us one motivation is to understand results as in \cite{KiFT}, where the philosophy is
to show that certain symplectic actions have to be Hamiltonian by using properties
of an associated 'generalized moment'.

\section{Convexity}
In this section the acting group is supposed to be a torus $S$, and we look at the
image of the generalized moment
$$\mu:M\to\C^*\times(S^1)^r\cong \C^*\times\R^r/\z^r$$
where $C$ is the maximal subtorus in $S$\ with Hamiltonian action, r its codimension in $S$.
Then, by the Atiyah-Guillemin-Sternberg convexity theorem, we know that $\mu_1(M)=\Delta$, a convex
polytope in $\C^*$.
\begin{theo}[Atiyah-Guillemin-Sternberg for $\mu$]
The image of $\mu$\ is\\ $\mu(M)=\Delta\times(S^1)^r$\ (convexity property),
and $\mu$\ is an open map from $M$\ to $\mu(M)$.

\end{theo}
{\bfseries{Tool :}} we will use the local form of the image of a (non-necessarily equivariant)
moment map for the $S$-action given by Theorem 32.3 of \cite{GuStSy} : the image of a
(non-necessarily invariant) neighbourhood $U$\ of a point $x\in M$\ for a local moment map
$\Phi$\ is given by
$$ \Phi(U)=U'\cap(\Phi(x)+s_1^\perp\oplus Co^+(\alpha_1,\ldots,\alpha_n)),$$
where $S_1$\ is the identity component of the the isotropy group at $x$, $s_1$\ its Lie
algebra and $s_1^\perp$\ its annihilator in $s^*$, $\alpha_1,\ldots,\alpha_n\in s_1^*\inc s^*$\ are the 
weights of the isotropy representation of $S_1$\ on $T_xM$, and \\$Co^+(\alpha_1,\ldots,\alpha_n)
\stackrel{def}{=}\{\sum_i\lambda_i\alpha_i|\lambda_i\in\R^+\}$.\\
Moreover, $\Phi$\ is open as a map $U\to \Phi(U)$\ (see \cite{HiNePl1994}).\\

{\bfseries{Proof of theo.4.1 :}} The inclusion $\mu(M)\inc\Delta\times(S^1)^r$\ is trivial.
We will show that a set $\mu(U)$\ described as $\Phi(U)$\ above is always open in
$\Delta\times(S^1)^r$, using again that McDuff moments have no local extrema. 
Combined with the fact that $\mu(M)$\ is closed in $\Delta\times(S^1)^r$, this will give the result.\\

Now if $\mu(U)$\ is not open, this means that $Co^+(\alpha_1,\ldots,\alpha_n)$\ is not $s_1^*$.
Let $V\inc s^*$\ be the subspace generated by $s_1^\perp\oplus Co^+(\alpha_1,\ldots,\alpha_n)$.
Then $W=V^\perp\inc s_1$\ is rationnaly generated (see a similar argument in the proof of lemma 2.2) so
it is the algebra of a subtorus $Z\inc S$. We claim that $Z\inc C$\ ; otherwise
there would exist some circle in $Z$, not contained in $C$.  This circle would have
a non-Hamiltonian action on $M$, whose McDuff moment would be constant on an open
set, a contradiction.  Now by looking at $\mu_1$\ and by usual Hamiltonian arguments,
one sees that $U$\ is the identity component of the generic stabilizer (corresponding to
the principal orbit type) of the points of $M$\ for the $C$-action.  This shows that
in fact $V=W^\perp$\ is independent of the point $x\in M$\ and the neighbourhood $U$\ chosen, and is given
by $>\Delta<\oplus \T^*$.  Let us write $S=W\oplus K$, with $K$\ the algebra of a complementary torus to $Z$\ in $S$.\\

We aim to show $\mu(U)=U'\cap(\Delta\oplus\T^*)=U'\cap V\cap(\Delta\oplus\T^*)$.  If
$V=s_1^\perp\oplus Co^+(\alpha_1,\ldots,\alpha_n)$, we are done.  Otherwise, by
the elementary properties of convex cones, $s_1^\perp\oplus Co^+(\alpha_1,\ldots,\alpha_n)$\ is
a finite intersection of half-spaces in $V$, of the form
$$H^+_y=\{v\in V|v(y)\geq 0\}\esp 0\neq y\in K\cap S_1$$
Again by rationality considerations, $>y<$\ is the algebra of a circle in $K$, and
if $y\not\in \C$\ we are led to a contradiction by creating an extremum for a 
McDuff moment.  So $y\in\C$, and the conclusion comes from the well-known openness of $\mu_1$.
\fin
\begin{cor}
If a n-torus $T^n$\ acts in a totally non-Hamiltonian way on a compact connected
symplectic manifold \sm, then $n\leq b_1(M)$.
\end{cor}
{\bfseries{proof :}}
Let $\mu :M\to (S^1)^r$\ be a generalized moment map.  Then by the surjectivity and 
the openness of $\mu$, it is possible to build n 1-cycles $\gamma_1,\ldots,\gamma_n$\ in $M$\ whose $\mu$-image are
(perhaps up to an integral multiple due to the non-connectedness of the fibers of
McDuff moments, which can force to make several revolutions before being back in the right component)
the 1-cycles in $T^n$\ corresponding to the n $S^1$-factors.  It is then easy to see that the
subgroup generated by $\gamma_1,\ldots,\gamma_n$\ in $H_1(M,\z)$\ is a free subgroup of rank n, giving $b_1(M)\geq n$.
\fin
\begin{rem}
This result has to be related to the "cohomologically free" actions in \cite{GiSACG}.
\end{rem}

From the cycle construction in Corollary 4.1, we also obtain the other
\begin{cor}
Let $\mu :M\to \C^*\times(S^1)^r$\ be a generalized moment for a symplectic torus action.
Let $(c,s)=(c,s_1,\ldots,s_{r-1})\in \C^*\times (S^1)^{r-1}$. Then there exists a
1-cycle $\gamma$\ in $M$\ with $\mu(\gamma)=(c,s_1,\ldots,s_{r-1})\times (S^1)$.
\end{cor}
This has an implication on the reduction side :
\begin{cor}
Suppose that $\mu$\ in Corollary 4.2 is naturally equivariant.  Decompose $\mu$\
as $\mu=(\mu_1,\mu_2)$\ where $\mu_1 :M\to \C^*\times(S^1)^{r-1}$\ and $\mu_2:M\to S^1$.
If $(c,s)$\ is a regular value of $\mu_1$\ and $M_{red}=\mu_1^{-1}(c,s)/{C\times (S^1)^{r-1}} $\
is a smooth reduced space, then the last $S^1$\ factor action induces a symplectic 
non-Hamiltonian action on $M_{red}$.
\end{cor}
{\bfseries{proof:}}
For clarity, we write $K$\ for the last $S^1$\ factor.\\
Let $i$\ denote the inclusion of $M_{cs}=\mu_1^{-1}(c,s)$\ in $M$, $\pi$\ the quotienting map
$\pi:M_{cs}\to M_{red}$. Then we see by standard arguments that $i^*(\mu_2)$\ is
a S-invariant function, furnishing a function $\tilde{\mu}_2:M_{red}\to S^1$\ which is a generalized moment
for the K-action on $M_{red}$.  The cycle $\gamma$\ in Corollary 4.2 sits in
$M_{cs}$, and projects in $M_{red}$\ to a cycle $\tilde{\gamma}$\ with 
$\tilde{\mu}_2(\tilde{\gamma})=S^1$, showing that really $\tilde{\mu}_2$\ is a McDuff moment 
corresponding to a non-Hamiltonian action.
\fin
\begin{rem}
Starting with a generalized moment $\mu:M\to\C^*\times(S^1)^r$, we can apply to 
each McDuff factor $\mu_2,\ldots,\mu_{r+1}:M\to S^1$\ the procedure of Remark 1.1,
to obtain a "moment" $\mu^{fc}:M\to\C^*\times(S^1)^r$\ where each $S^1$\ factor is fiber-connected.
It is not difficult to check that every result in this section remains true for $\mu^{fc}$.
\end{rem}

\begin{rem}
As reasonable extension of the result given here, we expect the Kirwan convexity theorem \cite{KiCPMM1984}
to be true in our context (see also \cite{GiCMM}).  Generalizations to symplectic orbifolds
are also imaginable.
\end{rem}
{\appendix
\section{Relation with a result of M.K. Kim}
As application of the "non-Hamiltonian heredity" in Corollary 4.3, we propose a reformulation
of the main theorem in \cite{KiFT}.  We have no familiarity with symplectic orbifolds, but we cannot
avoid their appearance in some arguments below, so we will make the following \\
{\bfseries{Assumption :}} the results used here and that we have proved in the manifold context are still
valid in the symplectic orbifold context.\\

The interesting main result stated by Kim is the
\begin{prop}
Let $\sigma:T^{(n)}\times M^{(2n+2)}\to M^{(2n+2)}$\ be an effective symplectic torus
action action on a compact connected symplectic manifold \sm, the supscripts indicating the respective dimensions.
Then $\sigma$\ is Hamiltonian whenever it possesses a fixed point.
\end{prop}

Suppose from now on that the proposition is false. Then a generalized moment
$\mu:M\to\C^*\times(S^1)^r$\ exists, and can be chosen naturally equivariant because
of the fixed point hypothesis. We will first prove a little lemma :
\begin{lem}
In the hypotheses of the proposition, there exists a decomposition $T=C\times (S^1)^{r-1}\times K\;\;,\;\;K\cong S^1$, and
a regular value $(c,s)\in \C^*\times(S^1)^{r-1}$ of $\mu_1$\ (notation of Corollary 4.3) such that $M_{cs}$\ contains a 
$C\times(S^1)^{r-1}$-orbit of $K$-fixed points. Equivalently, the
$K$-action on $M_{red}$\ has a fixed point.
\end{lem}
{\bfseries{Proof:}} We use the already invoked local form of any local moment map
near a fixed point $p\in M$, in a coordinate system $(x_i,y_i)$\ centered at $p$\ :
$$\mu(x_i,y_j)=\mu(0) + \demi\sum_{i=1}^{n+1} \alpha_i (x_i^2+y_i^2).$$
By usual arguments, the effectiveness of the action implies $>\{\alpha_i\}<=\T^*$.
Looking at the set $\{V_\alpha\}$\ of (n-1)-dimensional subspaces in $\T^*$\ of
the form\\ $>\alpha_1,\ldots,\hat{\alpha}_i,\ldots,\hat{\alpha}_j,\ldots,\alpha_{n+1}<$,
and their annihilators $>\xi_\alpha<\inc\T$, we see that the $\xi_\alpha$'s are circle generators
and at least one of thoses circles acts on $M$\ in a non-Hamiltonian way by dimensionality
reason. Pick one "non-Hamiltonian" $\xi_\alpha$, call $K$\ the corresponding circle. For notation ease,
reorder the coordinates so that the corresponding $V_\alpha$\ is $>\alpha_3,\ldots,\alpha_{n+1}<$. Remark that
every point in $M$\ having coordinates $(x_i,y_i)$\ with $x_1=x_2=y_1=y_2=0$\ is $K$-fixed. 
Let $F$\ be the relatively open cone $F=\{\sum_{i=3}^{n+1}\lambda^+_i\alpha_i|\lambda^+\in\R^+_0\}$\
Choose a subgroup $C\times (S^1)^{r-1}\inc T$\ such that $C\times (S^1)^{r-1}\times K$\ covers $T$,
and denote $i_1^*:\T^*\to \C^*\times(s_1^*)^{r-1}$\ the dual of the inclusion 
$\C\times(s_1)^{r-1}\inc \T$.  Let also $\mu=(\mu_1,\mu_2):M\to\C^*\times(S^1)^{r-1}\times K$\
be a corresponding naturally equivariant generalized moment.  Then, using for example
Sard's theorem to avoid singularities coming from outside our coordinate system, we
see that there exists a regular value $(c,s)$\ of $\mu_1$\ in the regio near $\mu_1(p)$\ corresponding
to $i_1^*(F)$, satisfying our requirements.
\fin
\begin{rem}
Again, this result remains valid if we replace $\mu$\ by a $\mu^{fc}$\ as in Remark 4.2.
\end{rem}  
{\bfseries{Proof of Kim's result :}} Still assuming that the result is false, 
consider a naturally equivariant "fiber-connected" moment $\mu^{fc}$.  Thanks to 
the results proved before, we see that a classical "reduction in stages" procedure 
allows to construct from $M$\ and $\mu^{fc}$\ a connected compact 4-dimensional orbifold
endowed with a non-Hamiltonian $K$-action possessing some fixed point.  Our assumption that
orbifolds behave like manifolds, applied to the McDuff result that symplectic $S^1$-actions
on compact connected 4-manifolds are Hamiltonian in the presence of fixed points (\cite{DuMMCA1988}),
furnishes our contradiction.
\fin
\begin{rem}
Without a way to produce connected reduced spaces, the arguments above fail
because there is a priori no reason for a "reduced" fixed point to be in a "non-Hamiltonian" reduced component.
This is why we used this somewhat inelegant trick with $\mu^{fc}$. 
\end{rem}
}

{\bfseries Acknowledgments :}\\
Whithout the insistant encouragements of my friends at U.L.B., this paper
would not exist. Many thanks to them all, with a special thought to
C\'eline Azizieh and Vincent Denolin.
\newpage
\nocite{GiCMH}
{\footnotesize{
\bibliographystyle{plain}

}}

\end{document}